\documentclass{article}
\usepackage[boxed]	{algorithm2e}
\usepackage		{amsmath}
\usepackage		{amsfonts}
\usepackage		{array}
\usepackage[english]	{babel}
\usepackage		{comment}
\usepackage		{cite}
\usepackage		{enumerate}
\usepackage[mathscr]	{eucal}
\usepackage		{graphicx}
\usepackage[utf8]	{inputenc}
\usepackage		{latexsym}
\usepackage		{longtable}
\usepackage		{multirow}
\usepackage		{rcs}
\usepackage		{theorem}
\usepackage		{xspace}

\newtheorem{theorem}{Theorem}
\newtheorem{definition}{Definition}
\newtheorem{lemma}{Lemma}

\newcommand{\Proof}{{\bf Proof}\xspace}

\DeclareMathOperator	{\cl}	{{\mathrm cl}}
\DeclareMathOperator	{\co}	{{\mathrm co}}

\DeclareMathOperator	{\intr}	{{\mathrm int}}

\newcommand{\R}{\mathbb{R}}

\newcommand{\cph}{CPH}

\newcommand{\MP}{MP\xspace}
\newcommand{\GMP}{GMP\xspace}
\newcommand{\mpdif}{\sim}
\newcommand{\nudif}{\div}
\newcommand{\cvE}{{\cal C}(E)}

\title{Yet Another Convex Sets Subtraction with Application in Nondifferentiable Optimization}
\author{Evgeni Nurminski \and Stan Uryasev
}
\date{\today}
\begin{document}

\maketitle

\begin{abstract}
This paper introduces a new subtraction operation for convex sets,
which defines their difference as a collection of inclusion-minimal
convex sets with appropriate definitions of linear operations on them.  
With these operations the set of collections becomes a linear vector space
with common zero and possibility to invert Minkowski summation.
As the demonstration of usability of this concept the Lipschitz continuity
of \(\epsilon\)-subdifferentials of convex analysis is proved in a novel way.
\end{abstract}

{\bf keywords:} convex sets, subtraction, differences, Minkowski-Pontryagin difference, linear space embedding

\section*{Introduction}
One of the most general approaches to solution of nonconvex nondifferentiable
optimization problems is the theory of quisidiffrentiable functions
whose directional derivatives at any point
of their domains can be represented
as the difference of two convex positively homogeneous (\cph)
functions of the direction.
This idea was originally put forward by B.N. Pshenichny
\cite{Pshenichny1980}
and later developed in many aspects by V.F.Demyanov et all
(see \cite{dem-vas-en-81,Demyanov_etal2014}).

As any \cph-function can be represented by its subdifferential at the origin
this representation produces a natural association of the directional
derivative with the pairs of convex sets
\([ \partial_\star f(x), \partial^\star f(x) ] = {\cal D} f(x)\)
which can be considered as upper and low subdiffrential parts of
quasidifferential \({\cal D} f(x)\) of \(f\) at a point \(x\).
It allows to formulate necessary optimality conditions
and with further refinements provide local optimization algorithms
by constructing descent directions \cite{Abbasov2013}.

Also this representation of directional derivatives alludes to
R{\aa}dstr{\"o}m 
\cite{Raadstroem1952}
embedding of the space of convex sets into the linear space of pairs of convex sets
factored by the equivalence relation 
for pairs \([A, B]\) of convex compacts:
\([A, B] \equiv [A+Z,B+Z]\) for any convex \(Z\)
and can be therefore be interpreted as the set difference of upper and lower subdifferentials.
Such interpretation is quite natural for nonconvex dc-optimization
\cite{HUrr-1985,Strekalovsky1990,Zhang2013}
and may bring some new ideas into optimization theory and algorithmica.

This work by R{\aa}dstr{\"o}m \cite{Raadstroem1952}
was used also by
Banks, Jacobs \cite{Banks+Jacobs1969} to develop a corresponding differential calculus,
Bradley, Datko
\cite{Bradley+Datko1977}
to study the differentiability properties of set-valued measures,
and
by Dempe, Pilecka
\cite{DempePileckaDemyanov}
in bilevel programming.
One should mention also the early work by Tjurin
\cite{Tjurin1965}
who made use of the embedding result
and studied the differentiability of set-valued maps given by systems of inequalities.
It was also extensively used in a series of works
\cite{Pallaschke+Urbanski2002,Scholtes1992,Pallaschke+Urbanski1993,Pallaschke+Urbanski1994,Pallaschke+Urbanski1996,Urbanski1996,Pallaschke+Urbanski2002}
to study algebraic structures associated with R{\aa}dstr{\"o}m embedding.

However embedding \(\cvE\) into the linear factor space of \(\cvE \times \cvE \) with the equivalence
relation
\( [X,Y] \mpdif [X + Z, Y + Z ] \)
for \(X, Y, Z\) in \(\cvE\) leaves open practical questions of checking this equivalence
for the particular pairs from \(\cvE \times \cvE \).
Possibly driven by this shortcomings a few other definitions of set difference were proposed.
Among the first alternative definitions of set difference one can mention
Hukuhara  \cite{Hukuhara}
who defined the difference of two sets \(X\) and \(Y\)
as a set \(Z\) such that \(X = Y + Z\).
This lead to a differential calculus very similar to the usual one for single-valued mappings
but the conditions for existence of such set \(Z\) are so restrictive that the concept is not widely applicable.
Slightly more general notion of set difference was suggested by
Pontryagin
\cite{Pontryagin1967}
and used in optimal control, see \cite{Kolmanovsky+Gilbert1998}
for the review of these investigations.

Another popular approach to define the difference of convex compacts is to relate it to the difference of
the support functions.
Than we are basically in the realm of positive homogeneous dc-functions
\cite{HUrr-1985}
and can use techniques from
functional analysis to define various metrics and develop calculus rules for
operations on differences.
In two related papers \cite{Baier2001,Baier2001_2}
one can find an extensive review of this direction together with constructive
rules for manipulation of convex and concave parts of dc-representations.

The specific feature of our approach is to define instead a difference of convex sets as a single object ---
a collection of convex sets and provide common linear operations for such object and demonstrate
how it can be used to invert widely used Minkowski summation.
For historical records it may be noticed that it was introduced in the technical report
\cite{Nurminski1982}.
The present paper refines this idea,
corrects certain inaccuracies and provides new proofs.
\section{Notations and Preliminaries}
\label{nots+pre}
The basic space that we shall consider is the finite-dimensional euclidean vector space $E$.
We are going to deal mainly with compact convex subsets of  $E$
with naturally defined operations of addition and multiplication by real numbers.
The space of such sets will be denoted as \(\cvE\).
Algebraically \(\cvE\) is the semigroup with respect to addition which is known
as Minkowski or vector sum: for \(X, Y \in \cvE\) the sum
\(X + Y = \{ x + y: x \in X, y \in Y\} \in \cvE \).
Multiplication is defined traditionally: for real-valued \(\alpha\) the set \(\alpha X = \{ \alpha x : x \in X \}\).
The space \(\cvE\), however, does not form a linear space due to the absence of an operation analogous to subtraction.
From many points of view it would be desirable to
 introduce such an operation as the inverse of the vector addition of sets,
as this is a well-defined operation which is used in many applications.
This would eventually amounts to embedding the space of convex sets into a linear space
and the most popular attempt of this kind has already been performed in  \cite{Raadstroem1952}.
Here we supply the theoretical embedding procedure with more structure and
describe some new applications.

Before discussion of subtraction of convex sets we should first explain
used notations and define related concepts.
We already defined in a common way the sum of sets and their multiplication by real numbers.
The inner product of two vectors \(x\) and \(y\) from \(E\) is denoted by \(xy\)
and the euclidean norm of the vector \(x\) is denoted as usual  by \(\|x\|=\sqrt{x x}\).
This notation is naturally extended for sets:
\( \|X\|=\sup_{x \in X \subset E}  \|x\|\)
with the triangle inequality
\(  \| X + Y  \| \leq   \| X \| + \| Y  \| \)
obviously holding for any \(X,Y \subset E\).

For some set \(G \subset E\) and vector \(y\), the expression  \(Gy\) denotes the set of inner products \(\{ g y, g \in G\}.\)

Special notation is used for certain sets : a singleton \(\{0\}\) is denoted by 0,
the closed unit ball is denoted by \(  B =\{  x : \|x\| \leq 1 \} \).

The Hausdorff distance between sets \(X\) and \(Y\) is defined as
\( d(X,Y) = \max \{d^0(X,Y),d^0(Y,X) \} \),
where
\( d^0(X,Y)=\max_{x \in  X} \min_{y \in  Y} \| x - y \|\).
The alternative definition 
\( d^0(X,Y)=\inf \{ \epsilon:  X \subset Y + \epsilon B,\,\epsilon \geq 0 \}\)
makes explicit use of Minkowski set addition.

The support function of a convex set \(X\) is denoted by:
\((X)_p= \sup_{x \in  X} px \)
and it is easy to check that \((X)_p\) is a \cph-function of \(p\).
Support functions have a number of other useful properties,
which relate geometrical features of convex sets to analytic properties of
\((X)_p \) .
One useful representation is
\(\|X\| = \sup_{p \in  B}  (X)_p \)
which suggests another definition of the Hausdorff distance:
\[ d(X,Y)=\sup_{p \in  B} (X)_p - \inf_{p \in  B}   (Y)_p\]
for \( Y \subset X\) which is known as Hormander equality.

Other notations used in the following sections include
\(\intr (X)\) --- the interior of set \(X\),
\(\co(X)\) --- the convex hull of set \(X\), and \(\cl(X)\) --- the closure of set \(X\).

We also make use of some well-known separation results for convex sets
(see f.i. classic \cite{Rockafellar1970})
which we formulate as follows.
\begin{theorem}
[separation theorem]
\label {tN1}
If \(X\) and  \(Y\)  are convex sets and \(\intr(X) \cap \intr(Y) \neq \emptyset\),
then there exists a vector \(p \neq 0 \) such that
\( (X)_p + (Y)_{-p} \leq 0 \).
\end{theorem}

However, a stronger result is more often used in practice.

\begin{theorem}
[strict separation theorem]
\label {tN2}
If  \(X\) and  \(Y\)  are convex sets such that \(\cl(X) \cap \cl(Y) = \emptyset \)
and at least one of these sets is bounded,
then there exists a vector \(p\) such that
\[(X)_p + (Y)_{-p} < 0.\]
\end{theorem}
This result is often used in the following form:
for a closed convex set \(X\) to contain zero vector \(0\)
it is necessary and sufficient that \((X)_p \geq 0\) for any \(p\). 

Next
we briefly review some important properties of the addition of convex sets
which will be useful in the discussion of subtraction that follows.

First, it is easy to check that the sum of two convex sets is also a convex set.
Also, for non-negative scalars \(\alpha\) and \(\beta\) 
\[\alpha X +  \beta  X = (\alpha+\beta)X\]
for any convex set \(X\)
This equality may fail if any of the scalars is negative or if set \(X\) is not convex.

Another nice feature of convex sets is that if \(X, Y, \) and \(Z\) are bounded convex
sets then  \(X \subset Y \) implies \( X + Z  \subset   Y + Z \) and
the converse is also true, i.e.,  \( X + Z  \subset   Y + Z \)  implies \(X \subset Y \).
\section{Generalized Minkowski-Pontryagin difference}
\label{gmp-diff}
One of the earliest definitions of set subtraction was given by
H. Hadwiger \cite{Hadwiger1950}
and called him ``Minkowski subtraction''
even if there is no evidence that Minkowski
actually defined it.
Later it was reintroduced by Pontryagin
\cite{Pontryagin1967}
and used in studies of optimal control problems,
so it makes sense to call it Minkowski-Pontryagin difference (\MP-difference for short).

The definition itself
goes like following:
\begin{definition}
[\MP-difference]
The \MP-difference of two convex sets \(X\) and \(Y\) is the set
\( X \mpdif Y = Z = \{ z: X \subset Y + z \}\).
\end{definition}

The \MP-difference of sets has a number of useful properties,
some of them are listed below:
\begin{itemize}
\item
The \MP-difference can be represented as the intersections (possibly empty) of translations:
\(X \mpdif Y = \cap_{x \in X} (Y - x) \)
\item
It can be used to approximate the minuend from above:
\(X \subset (X \mpdif Y) + Y \)
\item
\( 0 \in Y \) implies \( X \mpdif Y \subset X \)
\end{itemize}
The disadvantage of \(\mpdif\) is that \( (X \mpdif Y) + Y \neq X \)
so it does not invert summation of sets.
Also \MP-difference exists for rather narrow set of pairs \(X, Y\) where \(X\) is covered by
some translation of \(Y\).

In search for more general definition of set difference
we suggest another definition which is directly inspired by \MP-difference,
but compare to it is universal in a sense that it is applicable to any
pair of sets from \(\cvE\).

To do so we consider collections of convex closed bounded subsets of \(E\)
as a new object and define a set difference in the following manner .

\begin{definition}[Generalized \MP-difference]
For \(X, Y\) from \(\cvE\) the difference \(X \nudif Y\) is a collection of
inclusion-minimal \(Z \in \cvE\) such that \( X \subset Y + Z \).
\label{GMP-diff}
\end{definition}
The set of such collections forms the linear space by itself
with the natural definitions
of multiplication by real numbers as \( \alpha (X \nudif Y) = \alpha X \nudif \alpha Y \)
regardless of the sign of \(\alpha\),
and addition \( (X \nudif Y ) + ( Z \nudif W ) = (X + Z) \nudif (Y + W) \)
which will be discussed later on.

\begin{figure} 
\centering
\includegraphics[scale=.60]{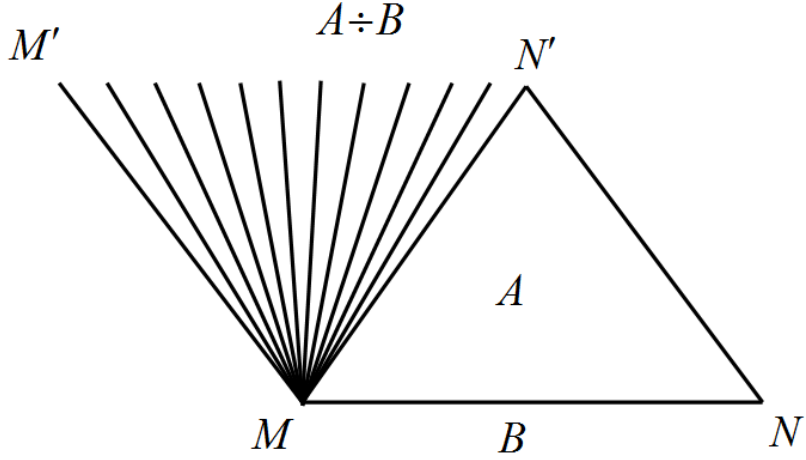} 
\caption{An example of the difference of two sets.} 
\label{fig:diference}
\end{figure}

The Definition
\ref{GMP-diff}
guarantees that \(X \nudif Y\) exists for any \(X, Y\)
as at least \(X \subset Y + \rho B\) for \(\rho\) large enough.
Also as any finite chain of \(Z_i, i = 1,2, \dots, k\) in \(\cvE\),
ordered by inclusion has a minimal element 
\( Z_{\min} = \cap_{i=1}^k Z_i \in \cvE\),
such minimal elements exist for the Definition \ref{GMP-diff}.
The properties of  \( X \nudif Y \)  are discussed in some detail in the next section;
here we just point out that Figure  \ref{fig:diference} provides the example of the case when
\(X \nudif Y\) contains more than one minimal element.

Here  set \(A\) is taken to be a triangle \(MN'N\)
( where the vertex \(M\) coincides with the origin )
and set \(B\) its base, the interval \(MN\).
The difference  \(A \nudif B\) 
is the collection of line segments connecting
the vertex \(M\) with any point on the interval \(M'N'\),
which is parallel to the base and has the same length.


It is easy to show that the difference  \(X \nudif Y\)
has also a number of important properties,
one of them is the ability to reproduce the difference of support functions
\((X)_p\) and \((Y)_p\) as it is
described in the following theorem.
\begin{theorem}
\label {tN3}
For any \(X, Y \in \cvE\)
\[ \inf_{Z \in X\nudif Y } (Z)_p  =    (X)_p -  (Y)_p \]
for any \(p\).
\end{theorem}
{\bf Proof.}
It is clear that at least \( (Y)_p + (Z)_p \geq (X)_p\) as for any \(Z \in X\nudif Y\)
\[   \inf_{Z \in X\nudif Y}  (Z)_p \equiv  \alpha \geq  (X)_p -  (Y)_p \equiv  \beta. \]
If this inequality is strict for some \(p'\) consider the intersection \(H_\rho\)
of the half-space \(H\) defined as follows :
\[ H = \left \{z:    p'z  \leq \frac {\alpha+\beta} {2} \right \}.  \]
with the ball \(\rho B\) with \(\rho\) arbitrary large:
\(H_\rho = \rho B \cap H\).
On one hand \( X \subset Y + H_\rho \)
because
\[  (Y + H_\rho)_p =   (Y)_p  + (H_\rho)_p  > (X)_p\]
for \( p \) not collinear with \(p'\) and \(\rho\) large enough.
At the same time
\[
\begin{array}{c}
(Y + H_\rho )_{p'} =   (Y)_{p'}  + (H_\rho)_{p'}=
(Y + H)_{p'} =   (Y)_{p'}  + (H)_{p'}=
\\
(X)_{p'} - \beta + \frac {\alpha+\beta} {2}=
(X)_{p'}  + \frac {\alpha-\beta} {2 } >(X)_{p'}.
\end{array}
\]
On other hand, there is no \( \bar Z \in X \nudif Y \) and \( \bar Z \subset H \).
If there were such a set, we would have
\[ \frac {\alpha+\beta}{2}=(H)_{p'} \geq (\bar Z)_{p'} \geq
\inf_{Z \in X \nudif Y} (Z)_{p'} = \alpha > \alpha- \frac {\alpha-\beta} {2 }= \frac {\alpha+\beta} {2 }
\]
and the theorem is proved by contradiction.

Next we describe some algebraic properties of subtraction and
relation of this operation to the Hausdorff distance between sets.
It is shown that subtraction of convex sets has quite standard algebraic properties,
although some of these ( for instance monotonicity ) are weaker than
the corresponding properties for real numbers.

It also demonstrates that the set of such collections forms the linear space
with the natural definitions
of multiplication by real numbers as \( \alpha (X \nudif Y) = \alpha X \nudif \alpha Y \)
regardless of the sign of \(\alpha\),
and addition \( (X \nudif Y ) + ( Z \nudif W ) = (X + Z) \nudif (Y + W) \).
As \(X = X \nudif 0 \) by definition we can sum  sets and collections:
\(X \nudif Y + Z = (X+Z) \nudif Y\).

If we define \(\inf_{Z \in X\nudif Y } (Z)_p \) to be
a support function of the collection \(X\nudif Y\) and denote it
as \((X\nudif Y)_p\) then this function becomes additive
as
\[
((X\nudif Y) + (Z \nudif W))_p = (X\nudif Y)_p + (Z \nudif W)_p 
\]
with respect to such definition of summation of collections.
The later makes it possible to relate our results to the approach of
Baier,Farkhi \cite{Baier2001,Baier2001_2}
who use decomposition of the set difference (in our understanding)
\(X \nudif Y\) onto convex and concave parts like
\(X \nudif Y = (X \nudif 0) + (0 \nudif Y)\) with summation defined in our paper
and incrementally constructing supports of these parts.

We can show useful cancellation and distributive properties of
the difference. 
\begin{lemma}  \label {LN1}
For sets \(X, Y, Z\) in \(\cvE\)
\[ X \nudif Y =(X + Z) \nudif (Y + Z)                       \]
\end{lemma}
{\bf Proof.}
This follows from the fact that for convex sets  \(X + Z \subset Y + Z  \) implies \( X  \subset Y \) and vice verse.
\begin{lemma}
[distributive law]
\label {LN2}
For \(X \in \cvE\) and
\( 0 \leq \gamma \leq 1 \)
\[ X \nudif \gamma X = (1 - \gamma) X. \]
\end{lemma}
{\bf Proof.}
This follows immediately from
\[ (X)_p  -  (\gamma X)_p = (X)_p - \gamma (X)_p = (1 -  \gamma) (X)_p = ((1-\gamma)X)_p,\]
where the last term is a \cph-function.
\begin{lemma}
[invertability]
\label {LN3}
\( 0= X \nudif Y \)
if and only if \( Y = X. \)
\end{lemma}
{\bf Proof.}
The equality \(  X \nudif X = 0 \) follows from Lemma \ref{LN2} with \(\gamma=1.\)
That this condition is sufficient can be proved in the following way:
\( 0 = X \nudif Y \)
immediately implies that
\( X  \subset  Y. \)

Furthermore, if there is a vector \(y\) such that \( y \in Y \) but  \( y \notin X  \),
then there exists a vector  \(p\) such that
\( (X)_p <  py \)
and
\[  0> (X)_p -  p y >   (X)_p -  (Y)_p = \inf_{Z \in  X \nudif Y } (Z)_p = (0)_p =0 \]
which proves the lemma by contradiction.
\begin{lemma} [monotonicity]
\label{LN4}
If  \( Y \subset X \)  then \(  0 \in X \nudif Y.  \)
\end{lemma}
{\bf Proof.}
Under the given conditions, for any \(Z \in X \nudif Y \)
\[ (Z)_p \geq  \inf_{Z \in X \nudif Y } (Z)_p =  (X)_p - (Y)_p \geq 0 \]
for any \(p\), and this proves the lemma.

\begin{figure} 
    \centering
    \includegraphics[scale=.60]{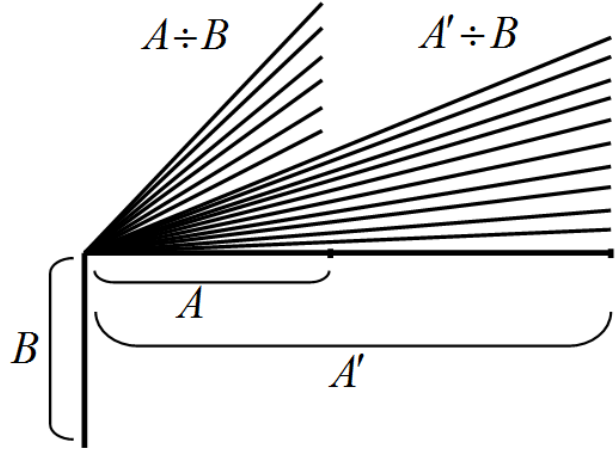} 
    \caption{Example showing that Lemma 4 cannot be generalized} 
   \label{fig:lemma4}
\end{figure}

Notice that the lemma in fact states that \(0 \in Z \) for any \(  Z \in   X \nudif Y\).
However, the counter-example given in Figure \ref{fig:lemma4} demonstrates that the generalization
of this lemma for \(  X \subset  X'  \) and differences \(X \nudif Y \) and \( X' \nudif Y  \) does not hold.

The norm of the set difference we define as
\[  \| X \nudif Y  \|   = \sup_{Z \in   X \nudif Y } \|Z\|.\]
Notice that for a convex set \( X \) there is a difference between
\( \| X \|=  \| X \nudif 0   \| \)
and the same set \(X\), considered as a collection of its elements.
It is easy to see that this collection can be represented as \( 0 \nudif X\) and so
it has the same norm  \( \| 0 \nudif X \| =\sup_{x \in X} \| x \|. \)
\section{Application of \GMP-difference to Convex Analysis}
\label{app}
The notion of the \GMP-difference of convex sets can help to study
 analytic properties of \(\epsilon\)-subdifferential mappings.
The concept of \( \epsilon$-subdifferential mapping proposed by
Rockafellar \cite{Rockafellar1970}
has proven to be very useful in convex nondifferentiable optimization.
This mapping is defined for any convex function \(f(x)\) and its value for a fixed point \( x\)
is a convex set \(\partial_{\epsilon}f(x )\) of vectors \( g\) such that
\[  \partial_{\epsilon}f(x ) \equiv \{g:   f(y) - f(x )   \geq g(y-x)  -\epsilon\}\]
for any \(y\) where \(\epsilon\) is a nonnegative constant.
It is worth to mention that \(\epsilon\)-subdifferential mappings were used
by Jeyakumar and Glover
\cite{Jeyakumar1996}
to derive global optimality conditions for nonconvex (dc) optimization.

There is a number of practical advantages in using \(\epsilon\)-subgradients in computational methods,
but the most interesting and promising feature of \(\epsilon\)-subgradiential mapping
is its richer analytic properties compared to the common subdifferential.
First of all it should be noted that for a fixed \(x\) the \(\epsilon\)-subdifferential mapping
as a multifunction of \(\epsilon \geq 0 \) has a convex graph
\(\{ (\epsilon, \partial_{\epsilon}f(x)), \epsilon \geq 0 \}\).
One of the of the earliest observations by
Asplund,Rockafellar
\cite{Asplund+Rockafellar1969}
was that
this mapping
unlike
classical subdifferential mapping
is continuous in Hausdorff metric.
Later results demonstrated that \(\epsilon\)-subdifferentia1 mapping 
is Lipschitz continuous in this metric as well for \(\epsilon > 0\)
and it was finally proved that it is even in some sense differentiable
\cite{Lemarechal+Nurminskii1980}.
It lead to the hope that the second derivatives will eventually be described
in a satisfactory way to provide foundation for Newton-like algorithms.
This hope did not realize so far so it make some sense to try different approaches.
Here we show again, using our definition of the subtraction of convex sets,
that \(\partial_{\epsilon}f(x )\)  is Lipschitz continuous for positive \(\epsilon\).
The result itself is known but the proof is simpler
and remarkably similar to the demonstration of the local Lipschitz property
of convex single-valued functions.

\begin{theorem}
[ Lipschitz Hausdorff continuity of \(\epsilon\)-subdifferential]
\label {tN6}
For a fixed \(x\) and $\epsilon>0$ for any \(\upsilon < \epsilon\) there exists \(L_{\epsilon, \upsilon}\) 
such that
\[
d(\partial_{\epsilon'}f(x ), \partial_{\epsilon''} f(x)) \leq L_{\epsilon, \upsilon} \vert \epsilon' - \epsilon'' \vert
\]
for any \(\epsilon', \epsilon''\) in \( [ \epsilon - \upsilon, \epsilon + \upsilon ]\).
\end{theorem}
\Proof 
To save on notation
we can will use the following shorthand
\( D( \epsilon)  \equiv \partial_{\epsilon}f(x ), D( 0)  \equiv \partial f(x ) \).
If  $0 < \epsilon - \upsilon < \epsilon' < \epsilon'' < \epsilon + \upsilon $.
then \( D( 0)  \subset    D(\epsilon')    \subset     D(\epsilon'') \)
and using convexity arguments for the graph of \(\partial_\epsilon f(x): \R_+  \to \cvE\) we have
\[
D(\epsilon')  \supset   \left  (1- \frac{\epsilon'' - \epsilon'}{\epsilon''}  \right )
D(\epsilon'') +   \frac{\epsilon''-\epsilon'}{\epsilon'} D(0).
\]
Adding $ \frac{\epsilon''-\epsilon'}{\epsilon''}    D(\epsilon'')$ to both sides yields
\[
D(\epsilon')  + \frac{\epsilon''-\epsilon'}{\epsilon''}  D(\epsilon'') \supset
D(\epsilon'')  +  \frac{\epsilon''-\epsilon'}{\epsilon''}  D(0)
\]
which can be rewritten as follows:
\[
\begin{array}{c}
D(\epsilon'')  +  \frac{\epsilon''-\epsilon'}{\epsilon''}  D(0) \subset
D(\epsilon')  + \frac{\epsilon''-\epsilon'}{\epsilon''}  D(\epsilon'')
\in
\\
D(\epsilon') + \frac{\epsilon''-\epsilon'}{\epsilon''}  D(0) +
\frac{\epsilon'' - \epsilon'}{\epsilon''} (D(\epsilon'')\nudif D(0)).
\end{array}
\]
Now we can drop $\frac{\epsilon''-\epsilon'}{\epsilon''}  D(0)  $ from both sides and obtain
\[
D(\epsilon'') \in  D(\epsilon')   +     \frac{\epsilon''-\epsilon'}{\epsilon''} (D(\epsilon'')\nudif D(0) ).
\]
By definition there is a set \( Z \in  D(\epsilon'')\nudif D( 0)$ such that
\( Z \subset (\|  D(\epsilon'')\nudif D(0) \|+\delta) B\)
for any \( \delta >0  $ and hence
\[
D(\epsilon'') \subset  D(\epsilon')   +     \frac{\epsilon''-\epsilon'}{\epsilon''} (\|D(\epsilon'')\nudif D(0)\|+\delta )B.
\]
The norm of the difference \( D(\epsilon'') \nudif D(0) \) is bounded from above
by a certain constant \( L_{\epsilon, \upsilon} \) which depends on \(\epsilon, \upsilon \)
so we can write
\[
\|D(\epsilon'') \nudif D(\epsilon')\| \leq L_{\epsilon, \upsilon} \frac{\epsilon''-\epsilon'}{\epsilon''}+\delta
\]
for arbitrary \( \delta>0 \),  which implies that
\[
\| D(\epsilon'')\nudif D(\epsilon') \| \leq L_{\epsilon, \upsilon} \frac{\epsilon''-\epsilon'}{\epsilon'}
\leq
\frac{L_{\epsilon, \upsilon}}{\epsilon-\upsilon}(\epsilon'' - \epsilon') \leq
L'_{\epsilon, \upsilon}(\epsilon'' - \epsilon')
\]
which demonstrates locally Lipschitz continuity of \(D(\epsilon)\)  in the Hausdorff metric
as a set-valued function of \(\epsilon\).

\textbf{Acknowledgments.}

For the first author this work is supported by the Ministry of Science and Education of Russian Federation,
project 1.7658.2017/6.7

\bibliographystyle{plain}
\bibliography{diff}
\end{document}